\documentclass[12pt]{article}%
\usepackage{amsfonts}
\usepackage{sw20bams}%
\usepackage{amsmath}%
\usepackage{amssymb}%
\usepackage{graphicx}
\providecommand{\U}[1]{\protect\rule{.1in}{.1in}}
\begin{document}

\title{Ptolemy Constants as Described by Eccentricity}
\author{Steven Finch}
\date{August 12, 2016}
\maketitle

\begin{abstract}
Let $J$ denote a simple closed curve in the plane. \ Let points $a$, $b$, $c$,
$d\in J$ occur in this order when traversing $J$ in a counterclockwise
direction. \ Define $p(a,b,c,d)$ to be the ratio of $ab\cdot cd+ad\cdot bc$ to
$ac\cdot bd$, where $zw$ denotes distance between $z$ and $w$. \ Define $P(J)$
to be the supremum of $p$ over all such points. \ Harmaala \&\ Kl\'{e}n
\cite{HK-Ptolemy} provided bounds on $P(J)$ when $J$ is an ellipse or
rectangle of eccentricity $\varepsilon$. \ We nonrigorously give formulas for
$P(J)$ here, in the hope that someone else can fill gaps in our reasoning.

\end{abstract}

\footnotetext{Copyright \copyright \ 2016 by Steven R. Finch. All rights
reserved.}Given parameters $0\leq\theta_{1}<$ $\theta_{2}<\theta_{3}%
<\theta_{4}<2\pi$ and $0\leq\varepsilon<1$, consider vertices
\[%
\begin{array}
[c]{ccc}%
v_{k}=\left(  \cos\theta_{k},\sqrt{1-\varepsilon^{2}}\,\sin\theta_{k}\right)
, &  & 1\leq k\leq4
\end{array}
\]
of a convex quadrilateral inscribed within the planar ellipse%
\[%
\begin{array}
[c]{ccc}%
x^{2}+\dfrac{y^{2}}{1-\varepsilon^{2}}=1 &  & \text{(with foci at }%
(\pm\varepsilon,0)\text{).}%
\end{array}
\]
The ratio%
\[
p=\frac{\left\vert v_{1}-v_{2}\right\vert \cdot\left\vert v_{3}-v_{4}%
\right\vert +\left\vert v_{1}-v_{4}\right\vert \cdot\left\vert v_{2}%
-v_{3}\right\vert }{\left\vert v_{1}-v_{3}\right\vert \cdot\left\vert
v_{2}-v_{4}\right\vert }
\]
involves lengths of sides in the numerator and lengths of diagonals in the
denominator. \ Let $P(\varepsilon)$ denote the supremum of the ratio over all
parameters $\theta_{1}$, $\theta_{2}$, $\theta_{3}$, $\theta_{4}$. \ It is
thought that $P(\varepsilon)$ measures the \textquotedblleft roundness of
planar curves\textquotedblright. \ Harmaala \&\ Kl\'{e}n \cite{HK-Ptolemy}
proved that%
\[
\frac{1}{2}\left(  \frac{1}{\sqrt{1-\varepsilon^{2}}}+\frac{\sqrt
{1-\varepsilon^{2}}}{1}\right)  \leq P(\varepsilon)\leq\csc\left(  \frac
{\pi\sqrt{1-\varepsilon^{2}}}{2}\right)
\]
but evidently did not tighten these bounds. \ 

Symbolic calculations of the gradient vector and Hessian matrix of $p$
indicate that%
\[
(\theta_{1},\theta_{2},\theta_{3},\theta_{4})=\left(  0,\frac{\pi}{2}%
,\pi,\frac{3\pi}{2}\right)
\]
corresponds to a local maximum of $p$, regardless of the value of
$\varepsilon$. \ For example, the Hessian matrix at this point is%
\[
\frac{-\varepsilon^{4}}{8\sqrt{1-\varepsilon^{2}}}\left(
\begin{array}
[c]{cccc}%
\dfrac{3-\varepsilon^{2}}{2-\varepsilon^{2}} & 0 & \dfrac{1-\varepsilon^{2}%
}{2-\varepsilon^{2}} & 0\\
0 & \dfrac{3-2\varepsilon^{2}}{\left(  2-\varepsilon^{2}\right)  \left(
1-\varepsilon^{2}\right)  } & 0 & \dfrac{1}{\left(  2-\varepsilon^{2}\right)
\left(  1-\varepsilon^{2}\right)  }\\
\dfrac{1-\varepsilon^{2}}{2-\varepsilon^{2}} & 0 & \dfrac{3-\varepsilon^{2}%
}{2-\varepsilon^{2}} & 0\\
0 & \dfrac{1}{\left(  2-\varepsilon^{2}\right)  \left(  1-\varepsilon
^{2}\right)  } & 0 & \dfrac{3-2\varepsilon^{2}}{\left(  2-\varepsilon
^{2}\right)  \left(  1-\varepsilon^{2}\right)  }%
\end{array}
\right)
\]
and all conditions of the multivariate second derivative test are clearly met.
\ Numerical optimization techniques suggest that this, in fact, corresponds to
a \textit{global} maximum. \ We do not see how to verify this rigorously. \ If
an analytical workaround could somehow be discovered, we would have%
\[
P(\varepsilon)=\dfrac{2-\varepsilon^{2}}{2\sqrt{1-\varepsilon^{2}}}
\]
for an ellipse of eccentricity $\varepsilon$, which is the lower bound given
in \cite{HK-Ptolemy}, Theorem 1.7.

Consider instead vertices $v_{1}$, $v_{2}$, $v_{3}$, $v_{4}$ of a convex
quadrilateral inscribed within the planar rectangle%
\[
\max\left\{  \left\vert x\right\vert ,\frac{\left\vert y\right\vert }%
{\sqrt{1-\varepsilon^{2}}}\right\}  =1.
\]
Cyclicity is assumed as before. \ This is analogous to the ellipse, although
the existence of sharp corners changes the nature of the analysis. \ Here we
have%
\[
P(\varepsilon)=\left\{
\begin{array}
[c]{lll}%
\sqrt{2} &  & \text{if }0\leq\varepsilon\leq\sqrt{3}/2,\\
\dfrac{\sqrt{1+4\left(  1-\varepsilon^{2}\right)  }}{2\sqrt{1-\varepsilon^{2}%
}} &  & \text{if }\sqrt{3}/2<\varepsilon<1
\end{array}
\right.
\]
for a rectangle of eccentricity $\varepsilon$, which again is the lower bound
given in \cite{HK-Ptolemy}, Corollary 4.8. \ The threshold $\varepsilon
=\sqrt{3}/2$ implies $\sqrt{1-\varepsilon^{2}}=1/2$, that is, a transition
occurs at a $2\times1$ rectangle.

The left-hand rectangle in Figure 1 shows an optimizing vertex configuration
for $0\leq\varepsilon\leq\sqrt{3}/2$; the right-hand rectangle shows an
optimizing vertex configuration for $\sqrt{3}/2\leq\varepsilon<1$. \ For the
former,%
\[%
\begin{array}
[c]{lll}%
v_{1}=\left(  1,\sqrt{1-\varepsilon^{2}}\right)  , &  & v_{2}=\left(
-1,\delta-\sqrt{1-\varepsilon^{2}}\right)  ,\\
v_{3}=\left(  -1,-\sqrt{1-\varepsilon^{2}}\right)  , &  & v_{4}=\left(
\delta-1,-\sqrt{1-\varepsilon^{2}}\right)
\end{array}
\]
give%
\[
\frac{\sqrt{4+\left(  \delta-2\sqrt{1-\varepsilon^{2}}\right)  ^{2}}%
\cdot\delta+\sqrt{\left(  \delta-2\right)  ^{2}+4\left(  1-\varepsilon
^{2}\right)  }\cdot\delta}{\sqrt{4+4\left(  1-\varepsilon^{2}\right)  }%
\cdot\sqrt{2}\,\delta}\rightarrow\sqrt{2}%
\]
as $\delta\rightarrow0^{+}$. \ For the latter,%
\[%
\begin{array}
[c]{lll}%
v_{1}=\left(  0,\sqrt{1-\varepsilon^{2}}\right)  , &  & v_{2}=\left(
-1,-\sqrt{1-\varepsilon^{2}}\right)  ,\\
v_{3}=\left(  0,-\sqrt{1-\varepsilon^{2}}\right)  , &  & v_{4}=\left(
1,-\sqrt{1-\varepsilon^{2}}\right)
\end{array}
\]
and the rest follows trivially. \ A sizeable variety of vertex configurations
need to be ruled out in order to verify global maximality.%
\begin{figure}[ptb]%
\centering
\includegraphics[
height=1.3439in,
width=5.8773in
]%
{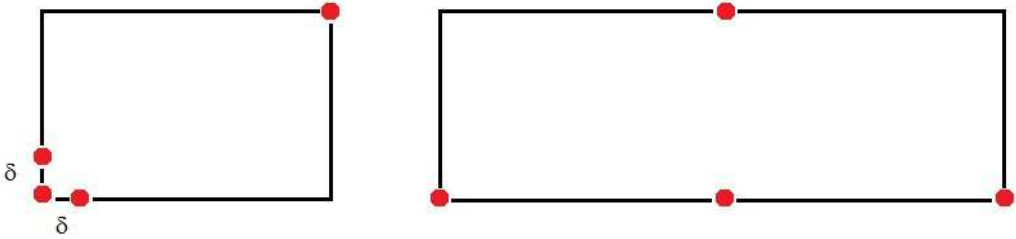}%
\caption{On the left are\ rectangles that are $2\times1$ or less eccentric.
\ On the right are rectangles that are $2\times1$ or more eccentric.}%
\end{figure}

Ptolemy constants remain open for a regular hexagon and for a Reuleaux
triangle, as well as for arbitrary convex quadrilaterals. \ Discovering these
could be a fruitful exercise in computer algebra.

\end{document}